\documentstyle[11pt]{article}
\newcommand{\ba}{\noindent $\begin{array}}
\newcommand{\ea}{\end{array}$}
\newcommand{\be}{\begin{equation}}
\newcommand{\ee}{\end{equation}}
\newcommand{\bd}{\begin{displaymath}}
\newcommand{\ed}{\end{displaymath}}
\newcommand{\beq}{\begin{eqnarray*}}
\newcommand{\eeq}{\end{eqnarray*}}
\newcommand{\beqn}{\begin{eqnarray}}
\newcommand{\eeqn}{\end{eqnarray}}

\begin{document}

\pagestyle{plain}

\begin{center}

{\large \bf COMPUTATIONAL EXPERIMENTS WITH  ABS ALGORITHMS 
\vskip2mm FOR KKT LINEAR SYSTEMS}
\footnote{Work partly supported by grant GA CR 201/00/0080
and by MURST 1997 and 1999 Programmi
di Cofinanziamento}

\end{center}

\begin{center}
{\bf E. Bodon
\footnote{Department of Mathematics, University of Bergamo, Bergamo 24129,
Italy (bodon@unibg.it)},
A. Del Popolo
\footnote{Department of Mathematics, University of Bergamo, Bergamo 24129,
Italy (delpopolo@unibg.it)},
L. Luk\v{s}an
\footnote{Institute of Computer Science, Academy of Sciences of the Czech
Republic,
Pod vod\'arenskou v\v e\v z\'\i\ 2, 182 07 Prague 8, Czech Republic
(luksan@uivt.cas.cz)} and
E. Spedicato
\footnote{Department of Mathematics, University of Bergamo, Bergamo 24129,
Italy (emilio@unibg.it)}}
\end{center}

\vspace{2mm}
\vspace{2mm}

\section{Introduction}

In this report, we present numerical experiments with three particular 
ABS algorithms:

\begin{itemize}

\item[(1)] The Huang or modified Huang algorithm,

\item[(2)] The implicit LU algorithm,

\end{itemize}

\noindent These algorithms have been used for 
finding a solution of the linear KKT system
$$
\left[ \begin{array}{cc} B & A^T \\ A & 0 \end{array} \right]
\left[ \begin{array}{c} x \\ y \end{array} \right] = 
\left[ \begin{array}{c} b \\ c \end{array} \right],
$$
where $A \in R^{m,n}$, $B \in R^{n,n}$ is symmetric, $b \in R^n$, 
$c \in R^m$, $x \in R^n$, $y \in R^m$ and $m \leq n$.

\noindent The considered algorithms belong to the basic ABS 
class of  algorithms for solving a linear system $Ax=b$, which is 
a special case of the scaled ABS class, given by the
following scheme:

\begin{itemize}

\item[(A)] Let $x_1 \in R^n$ be arbitrary and $H_1 \in R^{n,n}$ be  
nonsingular arbitrary. Set $i = 1$.

\item[(B)] Compute the residual $r_i = A x_i - b$. If $r_i = 0$, then 
stop, $x_i$ solves the problem. Otherwise compute $s_i = H_i A^T v_i$, 
where $v_i \in R^n$ is arbitrary, save that $v_1, \dots, v_i$ are linearly 
independent. If $s_i \neq 0$, then  go  to (C). If $s_i = 0$ and 
$r_i^T v_i = 0$, then set $x_{i+1} = x_i$, $H_{i+1} = H_i$ and go to (F). 
If $s_i = 0$ and $r_i^T v_i \neq 0$, then stop, the system is incompatible.

\item[(C)] Compute the search vector $p_i$ by
$$
p_i = H_i^T z_i,
$$
where $z_i \in R^n$ is arbitrary, save that $z_i^T H_i A^T v_i \neq 0$.

\item[(D)] Update the estimate of the solution by
$$
x_{i+1} = x_i - \alpha_i p_i,
$$
where the stepsize $\alpha_i$ is given by
$$
\alpha_i = r_i^T v_i / p_i^T A^T v_i. 
$$

\item[(E)] Update the Abaffian matrix by
$$
H_{i+1} = H_i - H_i A^T v_i w_i^T H_i / w_i^T H_i A^T v_i,
$$
where $w_i \in R^n$ is arbitrary, save that $w_i^T H_i A^T v_i \neq 0$.

\item[(F)] If $i = m$, then stop, $x_{i+1}$ solves the problem. Otherwise 
increment the index $i$ by one and go to (B).

\end{itemize}

 From (E), it follows by induction that $H_{i+1} A^T V_i = 0$ and  
$H_{i+1}^T W_i = 0$, where $V_i = [v_1, \dots, v_i]$ and 
$W_i = [w_1, \dots, w_i]$.
One can show that the implicit factorization $V_i^T AP_i = L_i$ holds, where
$P_i = [p_1, \dots, p_i]$ and $L_i$ is nonsingular lower triangular.
Moreover the general solution of the scaled subsystem $V_i^T A x =V_i^T b$ 
can be expressed in the form
\be
x = x_i + H_i^T q, 
\label{0}
\ee
where $q \in R^n$ is arbitrary (see \cite{ab2} for the proof).

The basic ABS class is the subclass of the scaled ABS class
obtained by taking $v_i = e_i$, $e_i$ being the $i$-th unitary vector 
($i$-th column of the unit matrix). In this case, residual $r_i = A x_i - b$ 
need not be computed in (B), $i$-th element $r_i^T e_i = a_i^T x_i - b_i$ 
suffices. 

This report is organized as follows. In Section 2, a short description of 
individual algorithms is given. Section 3 contains some details concerning
test matrices and numerical experiments. Section 4 discusses the results of the
numerical experiments.  The Appendix contains  all numerical results.
For a listing of the ABS source codes see [12]. For other numerical
results on ABS methods see [1,11].

\section{ABS algorithms for KKT equations}

In this report, we use  three linear ABS algorithms, two belonging to the
basic ABS class and one belonging to the scaled ABS class, for solving
the KKT linear systems exploiting their structure. To simplify
description, we will assume that $A \in R^{m,n}$, 
$m \leq n$, has full row rank so that $s_i \neq 0$ in (B). First we
describe the used linear solvers.

The Huang algorithm is obtained by the parameter choices $H_1 = I$, 
$v_i = e_i$, $z_i = a_i$, $w_i = a_i$. Therefore
\be
p_i = H_i a_i
\label{1}
\ee
and
\be
H_{i+1} = H_i - p_i p_i^T / a_i^T p_i,
\label{2}
\ee
From (2) and (3), it follows by induction that 
$p_i \in Range(A_i)$ and 
\be
H_{i+1} = I - P_i D_i^{-1} P_i^T. 
\label{3}
\ee
where $A_i = [a_1, \dots, a_i]$, $P_i = [p_1, \dots, p_i]$ and 
$D_i = diag(a_1^Tp_1, \dots, a_i^Tp_i)$. Moreover $H_i$ is symmetric,
positive semidefinite and idempotent (it is the orthogonal projection
matrix into $Null(A_{i-1})$). Since the requirement 
$p_i \in Null(A_{i-1})$ is crucial, we can improve orthogonality by 
iterative refinement $p_i = H_i^j a_i$, $j > 1$ (usually $j=2$), 
obtaining the modified Huang algorithm.  

The Huang algorithm can be used for finding the minimum-norm solution to 
the compatible underdetermined system $A x = b$, i.e. for minimizing 
$\| x \|$ s.t. $A x = b$. To see this, we use the Lagrangian function 
and convexity  of $\| x \|$. Then $x$ is a required solution if and only 
if $x = A^T u$  for some $u \in R^m$ or $x \in Range(A^T)$. But  
$$
x_{m+1} = x_1 - \sum_{i=1}^m \alpha_i p_i 
$$ 
by (D) and $p_i \in Range(A_i) \subset Range(A^T)$, so that if $x_1 = 0$, 
then $x_{m+1} \in Range(A^T)$. 

A short description of two versions of the Huang and modified Huang
algorithms follows.

\vspace{4mm}

\noindent {\bf Algorithm 1}

\noindent (Huang and modified Huang, formula (3).

\noindent Set $x_1 = 0$ and $H_1 = I$.

\noindent {\bf For} $i=1$ {\bf to} $m$ {\bf do}

\noindent Set $p_i = H_i a_i$ (Huang) or 

\noindent $p_i = H_i (H_i a_i)$ (modified Huang),

\noindent $d_i = a_i^T p_i$ and $x_{i+1} = x_i - ((a_i^T x_i - b_i)/d_i) p_i$.

\noindent If $i < m$, then set $H_{i+1} = H_i - p_i p_i^T / d_i$.

\noindent {\bf end do}

\vspace{4mm}

\noindent {\bf Algorithm 2}

\noindent (Huang and modified Huang, formula (4).

\noindent Set $x_1 = 0$ and $P_0$ empty.

\noindent {\bf For} $i=1$ {\bf to} $m$ {\bf do}

\noindent Set $p_i = (I - P_{i-1} D_{i-1}^{-1} P_{i-1}^T) a_i$ (Huang) or 

\noindent $p_i = (I - P_{i-1} D_{i-1}^{-1} P_{i-1}^T) 
		 (I - P_{i-1} D_{i-1}^{-1} P_{i-1}^T) a_i$  (modified Huang),
		 
\noindent $d_i = a_i^T p_i$ and $x_{i+1} = x_i - ((a_i^T x_i - b_i)/d_i) p_i$.

\noindent If $i < m$, then set $P_i = [ P_{i-1}, p_i ]$. 

\noindent {\bf end do}

\vspace{4mm}

The implicit LU algorithm is obtained by the parameter choices $H_1 = I$, 
$v_i = e_i$, $z_i = e_i$, $w_i = e_i$.  Since 
$W_i^T H_{i+1} = [I_i,0]^T H_{i+1} = 0$, the first $i$ rows of the 
Abaffian matrix $H_{i+1}$ must be zero. More precisely, the Abaffian matrix 
has the following structure, with $K_i \in R^{n-i,i}$ 
\be
H_{i+1} = \left[ \begin{array}{cc} 0 & 0 \\ K_i & I_{n-i} \end{array} \right].
\label{4}
\ee
Only the matrix $K_i$ has to be updated.        
Expression (5) implies that only the first $i$ components of the vector 
$p_i = H_i^T e_i$ can be nonzero and the $i$-th component is unity. Hence the 
matrix $P_i$ is unit upper triangular so that the implicit factorization 
$A = L P^{-1}$ is of the LU type with units on the diagonal.

A short description of two versions of the implicit LU algorithm follows (the 
second version is also called the implicit LX algorithm).

\vspace{4mm}

\noindent {\bf Algorithm 3}

\noindent (Implicit LU with explicit column interchanges).

\noindent Set $x_1 = 0$ and $H_1 = I$.

\noindent {\bf For} $i=1$ {\bf to} $m$ {\bf do}

\noindent Set $s_i = H_i a_i$ 

\noindent (only $(i-1)(n-i+1)$ nonzero elements of $H_i$ are used).

\noindent Determine $d_i = |s_{k_i}| = \max (|s_i^Te_j|, j=1, \dots ,n)$.

\noindent (only $n-i+1$ nonzero elements of $s_i$ are used). 

\noindent If $k_i \neq i$, then swap columns of $A$ and elements of 
	  $x$ and $s$ with these indices. 

\noindent Set $x_{i+1} = x_i - ((a_i^T x_i - b_i)/d_i) H_i^T e_i$ 

\noindent (only $i$ nonzero elements of $x_i$ are updated).

\noindent If $i < m$, then set $H_{i+1} = H_i - s_i e_i^T H_i^T / d_i$

\noindent (only $i(n-i)$ nonzero elements of $H_{i+1}$ are updated).

\noindent {\bf end do}

\vspace{4mm}

\noindent {\bf Algorithm 4}

\noindent (Implicit LX algorithm, or implicit LU with implicit column interchanges).

\noindent Set $x_1 = 0$, $H_1 = I$.

\noindent {\bf For} $i=1$ {\bf to} $m$ {\bf do}

\noindent Set $s_i = H_i a_i$ 
 
\noindent (only $(i-1)(n-i+1)$ nonzero elements of $H_i$ are used).

\noindent Determine $d_{k_i} = |s_{k_i}| = \max (|s_i^Te_j|, j=1, \dots ,n)$.

\noindent (only $n-i+1$ nonzero elements of $s_i$ are used). 

\noindent Set $x_{i+1} = x_i - ((a_i^T x_i - b_i)/d_{k_i}) H_i^T e_{k_i}$ 

\noindent (only $i$ nonzero elements of $x_i$ are updated).

\noindent If $i < m$, then set $H_{i+1} = H_i - s_i e_{k_i}^T H_i^T / d_{k_i}$

\noindent (only $i(n-i)$ nonzero elements of $H_{i+1}$ are updated).

\noindent {\bf end do}

\vskip4mm
\noindent We now consider how to use the above algorithm for solving
KKT systems, exploiting their structure. Consider the linear KKT system  
\be
B x + A^T y = b, 
\label{7}
\ee
\be
A x = c. 
\label{8}
\ee
The underdetermined equation (7) can be solved by an arbitrary
ABS method obtaining the particular solution $x_{m+1}$. Since
$H_{m+1} A^T = 0$, we get
\be
H_{m+1} B x = H_{m+1} b, 
\label{9}
\ee
by multiplying (6) by $H_{m+1}$. The coupled system (7)-(8) 
is overdetermined but compatible so that its unique solution can be 
found by an arbitrary ABS method. Since $rank(H_{m+1}) = n-m$, 
$m$ equations will be eliminated. Notice that solution of (7)-(8) 
is obtained in two steps. First we solve (7) to obtain 
$x_{m+1}$ and $H_{m+1}$. Then we construct $H_{m+1} B$ and $H_{m+1} b$ and 
starting with $x_{m+1}$ and $H_{m+1}$ we continue with the ABS method 
to solve (8). 

If we use the implicit LU algorithm, then matrix $H_{m+1}$ has the 
special form (3) and the first $m$ equations of (8) are trivially
satisfied, leaving a determined system consisting of (7) and the 
equation
\be
S_m B x = S_m b, 
\label{10}
\ee
where $S_m = [K_m, I_{n-m}]$. Moreover, substituting the general 
solution, see (1)
$$
x = x_{m+1} + H_{m+1}^T q, 
$$
with $q \in R^n$ arbitrary into (6) and using the special form of 
$H_{m+1}$, we obtain
\be
\left[ \begin{array}{cc} 0 & 0 \\ 0 & S_m B S_m^T \end{array} \right]
\left[ \begin{array}{c} q_1 \\ q_2 \end{array} \right] = 
\left[ \begin{array}{c} 0 \\ S_m (b - B x_{m+1}) \end{array} \right].
\label{11}
\ee
Thus $q_1$ can be chosen arbitrarily, 
e.g. $q_1 = 0$, and $q_2$ is a solution to the $n-m$ dimensional system 
\be
S_m B S_m^T q_2 = S_m (b - B x_{m+1}),
\label{12}
\ee
which can be solved by any ABS method. 

Once $x$ is determined, we obtain $y$ by solving the compatible
overdetermined system $A^T y = b - B x$. Since the equation $A x = c$
was solved beforehand, we are in the same position as in the case of 
the least-squares solution of an overdetermined linear system. Therefore, 
using $P_m$, we can construct the lower triangular matrix $L_m=AP_m$
and get the solution from the equation $L_m^T y = P_m^T (b - B x)$.

\section{Description of computational experiments}

Performance of ABS algorithms has been tested by using several types of 
ill-conditioned matrices. These matrices can be classified in the following
way. The first letter distinguishes matrices with integer 'I' and real 
'R' elements, both actually stored as reals in double precision arithmetic.      
The second letter denotes randomly generated matrices 'R' or matrices 
determined by an explicit formula 'D'. For randomly generated matrices,
a number specifying the interval for the random number generator follows,
while the name of matrices determined by the explicit formula contains
the formula number (F1, F2, F3). The last letter of the name denotes
a way for obtaining ill-conditioned matrices: 'R' - matrices with nearly
dependent rows, 'C' - matrices with nearly dependent columns, 'S' - nearly
singular symmetric matrices, 'B' -  both matrices in KKT system ill-conditioned.
More specifically:

\vspace{0.2cm}

\noindent \parbox{1.5cm}{IR500 } \parbox[t]{14.4cm}{Randomly generated matrices 
with integer elements uniformly distributed in the interval [-500,500].}

\noindent \parbox{1.5cm}{IR500R} \parbox[t]{14.4cm}{Randomly generated matrices 
with integer elements uniformly distributed in the interval [-500,500] perturbed 
in addition to have two rows nearly dependent.}

\noindent \parbox{1.5cm}{IR500C} \parbox[t]{14.4cm}{Randomly generated matrices 
with integer elements uniformly distributed in the interval [-500,500] perturbed 
in addition to have two columns nearly dependent.}

\noindent \parbox{1.5cm}{RR100 } \parbox[t]{14.4cm}{Randomly generated matrices 
with real elements uniformly distributed in the interval [-100,100]}.

\noindent \parbox{1.5cm}{IDF1  } \parbox[t]{14.4cm}{Matrices with elements 
$a_{ij} = |i-j|$, $1 \leq i \leq m$, $1 \leq j \leq n$ (Micchelli-Fiedler matrix).} 

\noindent \parbox{1.5cm}{IDF2  } \parbox[t]{14.4cm}{Matrices with elements 
$a_{ij} = |i-j|^2$, $1 \leq i \leq m$, $1 \leq j \leq n$.}

\noindent \parbox{1.5cm}{IDF3  } \parbox[t]{14.4cm}{Matrices with elements 
$a_{ij} = |i+j-(m+n)/2|$, $1, \leq i \leq m$, $1 \leq j \leq n$.}

\noindent \parbox{1.5cm}{IR50  } \parbox[t]{14.4cm}{Randomly generated matrices
with integer elements uniformly distributed in the interval [-50,50].}

\vspace{0.2cm}

Matrices with linearly dependent rows were obtained in the following way. 
The input data contain four integers which specify two row indices $i_1$, 
$i_2$, one column index $i_3$ and one exponent $i_4$. Then the row $a_{i_1}$
is copied into $a_{i_2}$. Furthermore $a_{i_1 i_3}$ is set to zero and
$a_{i_2 i_3}$ to $2^{-i_4}$. Similar procedures are used for columns 
and symmetric matrices.

Right hand sides were determined from the given solution vectors $x^{\star}$ 
by the formula $b = A x^{\star}$. Solution vectors were usually generated randomly 
with integer or real elements uniformly distributed in the interval [-10,10]. 

We have tested KKT systems with $m >> n/2$, $m = n/2$ and $m << n/2$. These systems 
were solved by using the modified Huang algorithm applied to the coupled system 
(7)-(8) and by two variants of the implicit LU algorithm with explicit column 
interchanges. The first variant is intended for 
solving the coupled system (7), 
(9) while the second one uses (7) together with the transformed system 
(11). For comparison, we have tested three additional methods implemented by 
using the LAPACK routines.

The first method, denoted by lu lapack, is in fact a direct solution of a complete
KKT system by using the Bunch-Parlett decomposition. This is carried-out by the 
LAPACK routine DSPSV. Notice that we deal with an $n+m$ dimensional indefinite
system in this case.   

Another approach, a range-space method, is based on the Bunch-Parlett decomposition
of the (possibly indefinite) matrix $B$. Using the LAPACK routine DSPSV, we obtain 
the solution of the matrix equation $B [\tilde{A}^T, \tilde{b}] = [A^T, b]$ together 
with the decomposition $B = L D L^T$, where $L$ is $n$-dimensional lower unit 
triangular and $D$ is $n$-dimensional block diagonal with the blocksize 1 or 2. 
Then the symmetric matrix 
$C = A \tilde{A}^T = A B^{-1} A^T$ is built and the solution $y$ to the equation 
$C y = \tilde{b} - c$ is found, again by using the LAPACK routine DSPSV. Finaly,
we solve the system $L D L^T x = b - A^T y$ by using the LAPACK routine DSPTRS.  

The last approach, a null-space method, is based on the RQ decomposition
\be
A = [0, R] Q, 
\label{13}
\ee
where $R$ is a $m$-dimensional upper triangular matrix and $Q$ is an $n$-dimensional 
orthogonal matrix. This RQ decomposition is obtained by using the LAPACK routine  
DQERQF. Premultiplying equation (6) by $Q$, we get $Q B Q^T Q x + Q A^T y = Q b$ 
or
\be
\left[ \begin{array}{cc} \tilde{B}_{11} & \tilde{B}_{12} \\ 
			 \tilde{B}_{21} & \tilde{B}_{22} \end{array} \right]
\left[ \begin{array}{c} \tilde{x}_1 \\ \tilde{x}_2 \end{array} \right] +
\left[ \begin{array}{c} 0 \\ R^T \end{array} \right] y = 
\left[ \begin{array}{c} \tilde{b}_1 \\ \tilde{b}_2 \end{array} \right].
\label{14}
\ee
where $\tilde{B} = Q B Q^T$, $\tilde{b} = Q b$ and $\tilde{x} = Q x$. Matrices $\tilde{B}$ 
and $\tilde{b}$ are obtained by using the LAPACK routine DORMRQ. Now $R \tilde{x}_2 = c$
by (7), (12) and 
$\tilde{B}_{11} \tilde{x}_1 = \tilde{b}_1 - \tilde{B}_{12} \tilde{x}_2$ 
by (13).
The last equation is solved by using the LAPACK routine DSYSV. Finally 
$R^T y = \tilde{b}_2 - \tilde{B}_{21} \tilde{x}_1 - \tilde{B}_{22} \tilde{x}_2$ and
$x = Q^T \tilde{x}$ can be obtained by the LAPACK routine DORMRQ.

 Notice that ABS algorithms were implemented in their basic form without
partitioning into block or other special adjustements serving for speed increase 
as done in the LAPACK software.

For each selected
problem, the type of matrix and the dimension is given. Furthermore, both the solution 
and the residual errors together with the detected rank and the computational time 
are given for each tested algorithm. Computational experiments were performed on a 
Digital Unix Workstation in the double precision arithmetic (machine epsilon equal
to about $10^{-16}$).

\vskip3mm
In the Appendix we give detailed results. The following  tables give
synthetic results for the 24 tested problems, the number at the intersection
of the $i$-th row with the $k$-th column indicating how many times the
algorithm at the heading of the $i$-th row gave a lower error 
than the algorithm at the heading of the $k$-th row (in case there is
a second number, this counts the number of cases when difference  was
less that one percent).
\newpage
\begin{verbatim}


       solution error  -  24  KKT linear systems

		  mod.   impl.  impl. lu      range  null
		  huang  lu8    lu9   lapack  space  space       total

    mod.huang            10     11     11     19     14          65
    impl.lu8      14             9/6   10     20     12          65/6
    impl.lu9      13      9/6           9     20     11          62/6
    lu lapack     13     14     15            22     15          79
    range space    5      4      4      2             4          19
    null space    10     12     13      9     20                 64



       residual error  -  24  KKT linear systems

		  mod.   impl.  impl. lu      range  null
		  huang  lu8    lu9   lapack  space  space       total

    mod.huang            15/1   16     11     20/1   6           68/2
    impl.lu8       8/1          15/2    8     19     5           55/3
    impl.lu9       8      7/2           7/2   20     6           48/4
    lu lapack     13     16     15/2          22     8/2         74/4
    range space    3/1    5      4      2            5           19/1
    null space    18     19     18     14/2   19                 88/2

\end{verbatim}
From the above tables and the results in the Appendix we can state
the following conclusions:

\begin{itemize}

\item[(1)] Implicit LU methods described in Subsection 2.3 are 
extremely suitable in term of total time cost for 
solving KKT systems with $m \sim n$. If $m << n$, then a direct solution of a complete
KKT system by using the Bunch-Parlett decomposition is faster. 

\item[(2)] The method based on equation (9) is faster if $m \sim n$ while the 
method based on equation (11) is more efficient if $m << n$.    

\item[(3)] The modified Huang method is usually very slow. But in the extremely 
ill-conditioned cases it gives the best accuracy. Moreover, the modified Huang method 
is able to detect the numerical rank correctly, which can lead to the substantial 
decrease of the computational time as can be observed from experiments with the matrix
IDF2. 

\item[(4)] In term of accuracy the range space method is the least accurate;
other methods have a comparable accuracy, with a marginal advantage for
lu lapack (except for the very ill conditioned problems, where mod.huang
gives the best results, up to four orders better).

\end{itemize}

\rm

\newpage

\section{Appendix: Test results for KKT linear systems}

\small

\begin{verbatim}

 matrix  dimension    method       solution  residual   rank     time
	   n    m                    error     error
 --------------------------------------------------------------------
 
 IR500   1000  900    mod.huang    0.50D-07  0.36D-14   1900   124.00
 IR500   1000  900    impl.lu8     0.70D-06  0.16D-13   1900    17.00
 IR500   1000  900    impl.lu9     0.71D-06  0.18D-13   1900    19.00
 IR500   1000  900    lu lapack    0.22D-05  0.96D-14   1900    56.00
 IR500   1000  900    range space  0.42D-01  0.52D-11   1900    79.00
 IR500   1000  900    null space   0.19D-06  0.20D-14   1900    85.00
condition number:     0.22D+10

 IR500   1200  600    mod.huang    0.32D-07  0.90D-14   1800   188.00
 IR500   1200  600    impl.lu8     0.83D-07  0.11D-13   1800    44.00
 IR500   1200  600    impl.lu9     0.89D-07  0.11D-13   1800    33.00
 IR500   1200  600    lu lapack    0.49D-07  0.41D-14   1800    47.00
 IR500   1200  600    range space  0.23D-02  0.55D-11   1800    64.00
 IR500   1200  600    null space   0.55D-08  0.26D-14   1800   105.00
condition number:     0.15D+09

 IR500   1500  200    mod.huang    0.15D-07  0.75D-14   1700   303.00
 IR500   1500  200    impl.lu8     0.76D-08  0.63D-14   1700    78.00
 IR500   1500  200    impl.lu9     0.92D-08  0.14D-13   1700    59.00
 IR500   1500  200    lu lapack    0.31D-08  0.91D-14   1700    40.00
 IR500   1500  200    range space  0.24D-05  0.16D-12   1700    41.00
 IR500   1500  200    null space   0.19D-08  0.25D-14   1700    78.00
condition number:     0.10D+08

 IR500R  1000  900    mod.huang    0.30D+01  0.49D-14   1900   124.00
 IR500R  1000  900    impl.lu8     0.80D-01  0.14D-13   1900    17.00
 IR500R  1000  900    impl.lu9     0.80D-01  0.10D-13   1900    19.00
 IR500R  1000  900    lu lapack    0.22D+01  0.10D-13   1900    55.00
 IR500R  1000  900    range space  0.74D+00  0.75D-11   1900    79.00
 IR500R  1000  900    null space   0.64D+01  0.28D-14   1900    85.00
condition number:     0.18D+17

 IR500R  1200  600    mod.huang    0.15D+00  0.81D-14   1800   187.00
 IR500R  1200  600    impl.lu8     0.35D-02  0.12D-13   1800    44.00
 IR500R  1200  600    impl.lu9     0.35D-02  0.14D-13   1800    34.00
 IR500R  1200  600    lu lapack    0.57D+00  0.64D-14   1800    47.00
 IR500R  1200  600    range space  0.16D+01  0.91D-11   1800    63.00
 IR500R  1200  600    null space   0.46D+00  0.20D-14   1800   105.00
condition number:     0.15D+16

 IR500R  1500  200    mod.huang    0.50D+00  0.98D-14   1700   302.00
 IR500R  1500  200    impl.lu8     0.81D-01  0.61D-14   1700    78.00
 IR500R  1500  200    impl.lu9     0.81D-01  0.83D-14   1700    59.00
 IR500R  1500  200    lu lapack    0.13D+01  0.16D-13   1700    39.00
 IR500R  1500  200    range space  0.27D+01  0.36D-12   1700    41.00
 IR500R  1500  200    null space   0.19D+00  0.31D-14   1700    78.00
condition number:     0.17D+16
\end{verbatim}           

\newpage

Test results for KKT linear systems - continued

\begin{verbatim}           
 matrix  dimension    method       solution  residual   rank     time
	   n    m                    error     error
 --------------------------------------------------------------------
 
 IR500S  1000  900    mod.huang    0.31D-06  0.36D-14   1900   125.00
 IR500S  1000  900    impl.lu8     0.41D-04  0.14D-13   1900    17.00
 IR500S  1000  900    impl.lu9     0.41D-04  0.30D-13   1900    21.00
 IR500S  1000  900    lu lapack    0.75D-05  0.19D-14   1900    56.00
 IR500S  1000  900    range space  0.34D+02  0.37D+01   1900    80.00
 IR500S  1000  900    null space   0.11D-04  0.19D-14   1900    89.00
condition number:     0.14D+12

 IR500S  1200  600    mod.huang    0.57D-08  0.37D-13   1800   189.00
 IR500S  1200  600    impl.lu8     0.12D-06  0.43D-12   1800    42.00
 IR500S  1200  600    impl.lu9     0.15D-06  0.16D-11   1800    34.00
 IR500S  1200  600    lu lapack    0.18D-08  0.19D-13   1800    49.00
 IR500S  1200  600    range space  0.18D+02  0.22D+01   1800    65.00
 IR500S  1200  600    null space   0.23D-07  0.19D-13   1800   110.00
condition number:     0.23D+09

 IR500S  1500  200    mod.huang    0.23D-08  0.76D-14   1700   304.00
 IR500S  1500  200    impl.lu8     0.13D-07  0.91D-14   1700    78.00
 IR500S  1500  200    impl.lu9     0.12D-07  0.14D-14   1700    60.00
 IR500S  1500  200    lu lapack    0.15D-08  0.11D-14   1700    41.00
 IR500S  1500  200    range space  0.14D+03  0.53D+02   1700    43.00
 IR500S  1500  200    null space   0.95D-08  0.19D-14   1700    77.00
condition number:     0.30D+08

 IR500B  1000  900    mod.huang    0.12D-01  0.83D-14   1900   125.00
 IR500B  1000  900    impl.lu8     0.14D-02  0.13D-13   1900    17.00
 IR500B  1000  900    impl.lu9     0.14D-02  0.14D-13   1900    20.00
 IR500B  1000  900    lu lapack    0.17D-02  0.45D-14   1900    58.00
 IR500B  1000  900    range space  0.26D+03  0.13D-02   1900    78.00
 IR500B  1000  900    null space   0.17D-01  0.23D-14   1900    88.00
condition number:     0.76D+14

 IR500B  1200  600    mod.huang    0.69D-03  0.60D-14   1800   189.00
 IR500B  1200  600    impl.lu8     0.15D-03  0.10D-13   1800    43.00
 IR500B  1200  600    impl.lu9     0.15D-03  0.11D-13   1800    33.00
 IR500B  1200  600    lu lapack    0.38D-02  0.20D-13   1800    47.00
 IR500B  1200  600    range space  0.66D+03  0.22D-03   1800    62.00
 IR500B  1200  600    null space   0.11D-01  0.21D-14   1800   108.00
condition number:     0.24D+14

 IR500B  1500  200    mod.huang    0.89D-03  0.11D-13   1700   309.00
 IR500B  1500  200    impl.lu8     0.23D-05  0.11D-13   1700    77.00
 IR500B  1500  200    impl.lu9     0.16D-05  0.40D-13   1700    57.00
 IR500B  1500  200    lu lapack    0.15D-03  0.42D-14   1700    39.00
 IR500B  1500  200    range space  0.22D+04  0.12D+00   1700    39.00
 IR500B  1500  200    null space   0.18D-02  0.34D-14   1700    77.00
condition number:     0.19D+13
\end{verbatim}           

\newpage

Test results for KKT linear systems - continued

\begin{verbatim}           
 matrix  dimension    method       solution  residual   rank     time
	   n    m                    error     error
 --------------------------------------------------------------------
 
 RR100   1000  900    mod.huang    0.57D-10  0.18D-12   1900   136.00
 RR100   1000  900    impl.lu8     0.64D-12  0.13D-13   1900    18.00
 RR100   1000  900    impl.lu9     0.50D-12  0.23D-13   1900    21.00
 RR100   1000  900    lu lapack    0.23D-11  0.11D-12   1900    60.00
 RR100   1000  900    range space  0.18D-08  0.10D-10   1900    87.00
 RR100   1000  900    null space   0.23D-12  0.26D-14   1900    91.00
condition number:     0.13D+04

 RR100   1200  600    mod.huang    0.92D-10  0.83D-12   1800   187.00
 RR100   1200  600    impl.lu8     0.13D-11  0.26D-13   1800    44.00
 RR100   1200  600    impl.lu9     0.16D-11  0.70D-13   1800    34.00
 RR100   1200  600    lu lapack    0.15D-11  0.11D-12   1800    47.00
 RR100   1200  600    range space  0.25D-09  0.19D-11   1800    64.00
 RR100   1200  600    null space   0.88D-13  0.24D-14   1800   105.00
condition number:     0.50D+03

 RR100   1500  200    mod.huang    0.23D-08  0.36D-10   1700   301.00
 RR100   1500  200    impl.lu8     0.10D-11  0.42D-13   1700    78.00
 RR100   1500  200    impl.lu9     0.33D-11  0.10D-12   1700    59.00
 RR100   1500  200    lu lapack    0.62D-12  0.11D-12   1700    39.00
 RR100   1500  200    range space  0.28D-10  0.87D-12   1700    41.00
 RR100   1500  200    null space   0.96D-13  0.23D-14   1700    78.00
condition number:     0.15D+03

 IDF1    1000  900    mod.huang    0.85D-10  0.18D-14   1900   133.00
 IDF1    1000  900    impl.lu8     0.38D-10  0.13D-14   1900    17.00
 IDF1    1000  900    impl.lu9     0.36D-10  0.14D-14   1900    20.00
 IDF1    1000  900    lu lapack    0.48D-10  0.35D-14   1900    59.00
 IDF1    1000  900    range space  0.79D-10  0.18D-14   1900    85.00
 IDF1    1000  900    null space   0.66D-10  0.21D-14   1900    92.00
condition number:     0.80D+04

 IDF1    1200  600    mod.huang    0.16D-08  0.27D-13   1800   182.00
 IDF1    1200  600    impl.lu8     0.55D-10  0.15D-14   1800    43.00
 IDF1    1200  600    impl.lu9     0.49D-10  0.15D-14   1800    33.00
 IDF1    1200  600    lu lapack    0.28D-10  0.19D-14   1800    47.00
 IDF1    1200  600    range space  0.87D-10  0.16D-14   1800    63.00
 IDF1    1200  600    null space   0.64D-10  0.25D-14   1800   105.00
condition number:     0.96D+04
    
 IDF1    1500  200    mod.huang    0.33D-07  0.42D-12   1700   303.00
 IDF1    1500  200    impl.lu8     0.52D-10  0.19D-14   1700    77.00
 IDF1    1500  200    impl.lu9     0.63D-10  0.15D-14   1700    58.00
 IDF1    1500  200    lu lapack    0.28D-10  0.15D-14   1700    39.00
 IDF1    1500  200    range space  0.11D-09  0.18D-14   1700    41.00
 IDF1    1500  200    null space   0.15D-09  0.33D-14   1700    78.00
condition number:     0.12D+05
\end{verbatim}           

\newpage

Test results for KKT linear systems - continued

\begin{verbatim}           
 matrix  dimension    method       solution  residual   rank     time
	   n    m                    error     error
 --------------------------------------------------------------------
 
 IDF2    1000  900    mod.huang    0.55D+01  0.23D-14     16    24.00
 IDF2    1000  900    impl.lu8     0.44D+13  0.21D-03   1900    18.00
 IDF2    1000  900    impl.lu9     0.12D+15  0.80D-02   1900    21.00
 IDF2    1000  900    lu lapack    0.25D+03  0.31D-13   1900    62.00
 IDF2    1000  900    range space  0.16D+05  0.14D-11   1900    87.00
 IDF2    1000  900    null space   0.89D+03  0.15D-12   1900    93.00
condition number:     0.26D+21

 IDF2    1200  600    mod.huang    0.62D+01  0.20D-14     17    36.00
 IDF2    1200  600    impl.lu8     0.22D+07  0.10D-08   1800    44.00
 IDF2    1200  600    impl.lu9     0.21D+06  0.56D-09   1800    33.00
 IDF2    1200  600    lu lapack    0.10D+03  0.79D-14   1800    47.00
 IDF2    1200  600    range space  0.11D+05  0.15D-11   1800    63.00
 IDF2    1200  600    null space   0.38D+04  0.13D-12   1800   105.00
condition number:     0.70D+20

 IDF2    1500  200    mod.huang    0.38D+01  0.22D-14     17    68.00
 IDF2    1500  200    impl.lu8     0.31D+04  0.76D-11   1700    79.00
 IDF2    1500  200    impl.lu9     0.58D+03  0.57D-12   1700    59.00
 IDF2    1500  200    lu lapack    0.28D+03  0.19D-13   1700    39.00
 IDF2    1500  200    range space  0.99D+05  0.16D-10   1700    41.00
 IDF2    1500  200    null space   0.25D+03  0.36D-13   1700    78.00
condition number:     0.28D+20

 IR50    1000  900    mod.huang    0.76D+00  0.48D-14   1702   143.00
 IR50    1000  900    impl.lu8     0.29D+16  0.28D+00   1900    18.00
 IR50    1000  900    impl.lu9     0.13D+16  0.15D+00   1900    21.00
 IR50    1000  900    lu lapack    0.10D+04  0.59D-13   1900    60.00
 IR50    1000  900    range space  0.24D+03  0.17D-10   1900    84.00
 IR50    1000  900    null space   0.18D+18  0.11D+02   1900    91.00
condition number:     0.27D+19

 IR50    1200  600    mod.huang    0.30D+00  0.84D-14   1751   186.00
 IR50    1200  600    impl.lu8     0.44D+14  0.53D-02   1800    44.00
 IR50    1200  600    impl.lu9     0.23D+14  0.33D-02   1800    34.00
 IR50    1200  600    lu lapack    0.63D+02  0.15D-13   1800    47.00
 IR50    1200  600    range space  0.10D+03  0.28D-09   1800    63.00
 IR50    1200  600    null space   0.82D+15  0.66D-01   1800   105.00
condition number:     0.14D+18

 IR50    1500  200    mod.huang    0.33D-07  0.68D-14   1700   302.00
 IR50    1500  200    impl.lu8     0.73D-08  0.76D-14   1700    78.00
 IR50    1500  200    impl.lu9     0.63D-07  0.99D-14   1700    59.00
 IR50    1500  200    lu lapack    0.58D-09  0.63D-14   1700    39.00
 IR50    1500  200    range space  0.20D-05  0.36D-13   1700    41.00
 IR50    1500  200    null space   0.11D-07  0.36D-14   1700    78.00
condition number:     0.12D+08

\end{verbatim}           

\rm
\end{document}